\documentclass[11pt,a4paper]{amsart}
\usepackage{amssymb,amscd}
\usepackage[dvipdfm]{graphicx}
\usepackage[dvipdfm,pageanchor,
bookmarks=true,bookmarksnumbered=true,bookmarkstype=toc,%
colorlinks=true,linkcolor=red,citecolor=blue,urlcolor={0 0 1},%
backref,bookmarks,breaklinks,linktocpage,%
pdftitle={High-codimensional knots spun about manifolds},%
pdfsubject={DRAFT},%
pdfauthor={Dennis Roseman and Masamichi Takase},%
pdfkeywords={spinning, Haefliger knot, long knot}%
]{hyperref}
\hypersetup{pdfstartview={FitH -32768}}
\DeclareGraphicsExtensions{.eps}
\usepackage{times,mathptmx}
\usepackage[psamsfonts]{eucal}

\newtheorem{thm}{Theorem}[section]

\newtheorem{prop}[thm]{Proposition}
\theoremstyle{definition}
\newtheorem{definition}[thm]{Definition}

\newtheorem{remark}[thm]{Remark}

\newcommand{\Z}{\mathbb{Z}}
\newcommand{\R}{\mathbb{R}}
\newcommand{\T}{\mathbb{T}}
\newcommand{\co}{\colon\thinspace}
\newcommand{\puncT}{\stackrel{\circ}{T}}

\newcommand{\K}{\mathcal{K}}
\newcommand{\gr}{\mathrm{gr}}

\begin{document}
\title{High-codimensional knots spun about manifolds}
\author{Dennis Roseman}
\address{%
Department of Mathematics, University of Iowa, 14 MacLean Hall,
Iowa City, 52242-1419 U.S.A.}
\urladdr{http://www.math.uiowa.edu/{\textasciitilde}roseman/}
\email{roseman@math.uiowa.edu}
\author{Masamichi Takase}
\address{%
Department of Mathematical Sciences, 
Faculty of Science, Shinshu University, 
Matsumoto, 390-8621 JAPAN}
\email{takase@math.shinshu-u.ac.jp}
\thanks{The second-named author is partially supported
by the Iwanami Fujukai Foundation.}

\pagestyle{headings}

\begin{abstract}
Using spinning we analyze in a geometric way Haefliger's smoothly 
knotted $(4k-1)$-spheres in the $6k$-sphere. 
Consider the $2$-torus standardly embedded in the $3$-sphere, which
is further standardly embedded in the $6$-sphere. At each point of
the $2$-torus we have the normal disk pair: a $4$-dimensional disk and
a $1$-dimensional proper sub-disk. We consider an isotopy
(deformation)
of the normal $1$-disk inside the normal $4$-disk, by using a map
from the $2$-torus to the space of long knots in 4-space, first
considered by Budney. We use this isotopy in a construction called
\textit{spinning about a submanifold} introduced by the first-named
author. Our main observation is that the resultant spun knot
provides a generator of the Haefliger knot group of knotted
$3$-spheres in the $6$-sphere.
Our argument uses an explicit construction of
a Seifert surface for the spun knot and works also for
higher-dimensional Haefliger knots.
\end{abstract}
\maketitle

\section{Introduction}\label{sect:introduction}
Various kinds of spinning constructions, all of which stem from
Artin's original construction \cite{AR}, are now basic tools in the
study of high-dimensional knots in \textit{codimension two}, that
is, in the study of embeddings of $n$-manifolds in the
$(n+2)$-sphere. 

On another front,  Haefliger found smooth knots in
\textit{codimensions greater than two}  \cite{HAE1,HAEF2}. He
showed that the group $C_n^q$ of isotopy classes of smooth
embeddings of the $n$-sphere $ S^n$ into $(n+q)$-sphere $ S^{n+q}$
is often non-trivial even when $q\ge3$.
Although many spinning constructions can be applied
for such ``high-codimensional'' knots, there have been very few
related studies \cite{budney,h-s,milgram}.

Budney \cite{budney} gave a new description of a generator of
the Haefliger knot group $C^3_3$.
This was related to  his study on the space $\mathcal{K}_{4,1}$ of
``long'' knots in $4$-space. Here, $\mathcal{K}_{n,j}$ is the space
of \textit{long knots} --- smooth embeddings of
$\R^j\hookrightarrow\R^n$ which are standard outside the unit disk
in $\R^j$.
We note that the restriction of a long embedding to the
unit disk gives a properly embedded $j$-disk in the unit $n$-disk 
--- our construction will refer to it. 
In his paper, a map (the resolution map) $\Phi\co T^2\to\K_{4,1}$ from the $2$-torus
$T^2$, defined in an ingenious way (see \S\ref{subsect:budney}),
plays a key role. This map gives rise to generators of certain
homotopy groups
and of the Haefliger group,
via successive graphing constructions.
The geometric nature of graphing constructions here
is a high-codimensional version of
Litherland's deform-spinning.

We give yet another description of a
generator of
the Haefliger knot group $C^3_3$ 
in terms of the notion \textit{spinning about a submanifold},
introduced by the first-named author \cite{DR-spin}.
Our main result is as follows (see \S\ref{sect:main} for the
details). Consider the $2$-torus $T^2$ standardly embedded in $
S^3$, which is further standardly embedded in $ S^6$. At each point
of $T^2$, we have the normal $4$-disk $ D^4$ to $T^2\subseteq S^6$
and the normal $1$-disk $ D^1$ to $T^2\subseteq S^3$. Then, we
``deform'' $ D^1$ inside $ D^4$, by using $\Phi$. Namely, at each
point of $\tau\in T^2$, we replace the standard disk pair $( D^4,
D^1)$ with the new disk pair $( D^4,\Phi(\tau)( D^1))$.
We then show that  the resultant smoothly embedded $3$-sphere
$\Sigma_\Phi^3\subseteq S^6$ represents a generator of the Haefliger
knot group $C_3^3$.

Our study is motivated by Budney's construction and in particular
his resolution map $\Phi\co T^2\to\K_{4,1}$. Our approach is
different. We use a basic technique in the spirit of codimension two
knot theory --- examination of how the homology of a Seifert surface
relates the knot complement. It is very geometric, uses higher-dimensional
visualization, does not involve any homotopy groups and might
be useful for more general high-codimensional knots.

Additionally, all of Budney's and our arguments work for
higher-dimensional Haefliger knot groups $C^{2k+1}_{4k-1}$, $k\ge2$.
We need just to consider the triple
\[ S^{2k-1}\times S^{2k-1}\subseteq S^{4k-1}\subseteq S^{6k}\]
and use everywhere a higher-dimensional Budney map
\cite[\S5]{budney}
\[ S^{2k-1}\times S^{2k-1}\to\mathcal{K}_{2k+2,1}\]
instead of $\Phi\co T^2= S^1\times S^1\to\K_{4,1}$.

Throughout the paper, we work in the smooth category; all manifolds
and mappings are supposed to be differentiable of class $C^\infty$,
unless otherwise stated. We use the symbol `$\cong$' for a group
isomorphism and `$\approx$' for a diffeomorphism.  Some graphics
are shown piecewise smooth --- they correspond to unique
smooth manifolds by rounding of corners.
%

The second-named author would like to thank
Ryan Budney for fruitful conversations.

\section{Preliminaries}
\subsection{Haefliger's knots}
Haefliger showed in \cite{HAE1,HAEF2} that the group $C^n_q$ of
smooth isotopy classes of smooth embeddings of the $n$-sphere $ S^n$
in the $(n+q)$-sphere $ S^{n+q}$ is often non-trivial even when the
codimension $q$ is greater than $2$. In a particular case, for each
$k\ge1$ the group $C^{4k-1}_{2k+1}$ of smooth isotopy classes of
smooth embeddings $ S^{4k-1}\hookrightarrow S^{6k}$ forms the
infinite cyclic group $\Z$. This is in contrast with Zeeman's
unknotting theorem \cite{z1} claiming that any $n$-sphere is
unknotted in the \emph{piecewise linear} sense in the $(n+q)$-sphere
if $q>2$.

According to \cite{boechat,b-h,mT1,mT2}, the smooth isotopy class of
Haefliger's knot
can be read off from geometric characteristics of
its Seifert surface.
When $k=1$,
we have the following (\cite{g-m,b-h}).

\begin{thm}\label{thm:invariant}
Every embedding $F\co S^3\hookrightarrow S^6$ has a Seifert
surface $\widetilde{F}\co V^4\hookrightarrow S^6$ and
\begin{eqnarray*}
\Omega(F)
&=&-\frac{1}{8}(\sigma(V^4)-e_{\widetilde{F}}\smile e_{\widetilde{F}})
\end{eqnarray*}
gives the isomorphism $\Omega\co C^3_3\to\Z$,
where $e_{\widetilde{F}}$ denotes the normal Euler class of $\widetilde{F}$.
\end{thm}

\subsection{Budney's isomorphism}\label{subsect:budney}
Several papers including \cite{boechat,b-h,g-m,mT1,mT2,skopenkov}
give geometric descriptions for the Haefliger knot groups
$C^{2k+1}_{4k-1}$ and their generators. Recently, Budney
\cite{budney} has given a description of such Haefliger knots in
terms of the space of long knots, which is related to the
Litherland-type deform-spinning.

Let $\K_{n,j}$ be the space of smooth embeddings (\textit{long
embeddings}) $\R^j\hookrightarrow\R^n$ being the standard inclusions on
$|x|\ge1$ for $x\in\R^j$. Then, Budney \cite{budney}, using
results in \cite{goodwillie}, showed that there is an isomorphism
$\pi_2 \K_{4,1} \to C^3_3$.
His isomorphism gives a new description of 
the Haefliger group $C^3_3\approx\Z$.

We briefly review Budney's construction. Note that he deals with
more general cases in \cite[\S5]{budney}.

\subsubsection{The graphing map}
Consider the ``graphing'' map $\gr_1\co\Omega\K_{n-1,j-1}\to\K_{n,j}$
defined by
\[(\gr_1 f)(t_0,t_1,\cdots,t_{j-1})
= \left(t_0,f(t_0)\left(t_1,\cdots,t_{j-1}\right)\right).\]
Budney \cite{budney} showed that the maps
$\Omega^2\K_{4,1}\to\Omega\K_{5,2}\to\K_{6,3}$
between the loop spaces induce isomorphisms
\[
\pi_2\K_{4,1}\cong\pi_1\K_{5,2}\cong\pi_0\K_{6,3}.
\]
Furthermore, the group $\pi_0\K_{6,3}$ is isomorphic to $C^3_3$
due to \cite{HAEF2}.

\subsubsection{The resolution map $\Phi\co T^2\to\K_{4,1}$}\label{subsubsect:resolution} 
\begin{figure}[tbhp]
\begin{center}
\includegraphics[width=.8\textwidth,keepaspectratio]{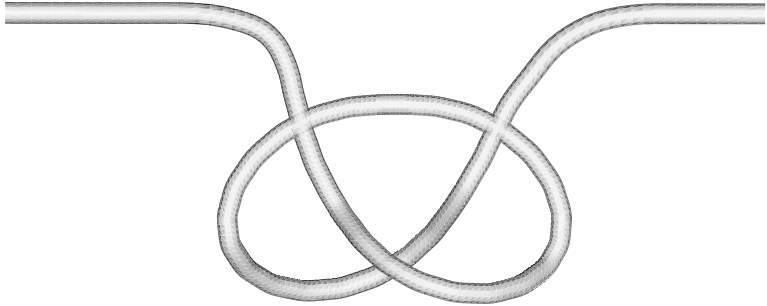}
\end{center}
\caption{An immersed long line with two intersections, thickened and shown as a thin tube.}\label{fig:resolution}
\end{figure}

Take a ``long'' immersion $f\co\R\to\R^3 \subseteq\R^4$
such that $f(t)=(t,0,0,0)$ for $|t|>1$ and has two
double points $f(t_1)=f(t_3)$, $f(t_2)=f(t_4)$ with
$-1 < t_1 < t_2 < t_3 < t_4 < 1$ and such that
$df(T_{t_1}\R)\cap df(T_{t_3}\R)=df(T_{t_2}\R)\cap df(T_{t_4}\R)=\{0\}$
(see Figure~\ref{fig:resolution}).
At $f(t_1)=f(t_3)$ and $f(t_2)=f(t_4)$,
we have the $2$-dimensional normal complements $P_1$
to $df(T_{t_1}\R)\oplus df(T_{t_3}\R)$
and $P_2$ to $df(T_{t_2}\R) \oplus df(T_{t_4}\R)$,
respectively.

Let $ S_1$ and $ S_2$ be the unit $1$-dimensional sphere in $P_1$
and $P_2$ respectively.
Given $(\theta,\psi)\in S_1\times S_2$, we perturb a small
neighborhood of $f(t_1)$ in direction $\theta$ using a bump function
and a small neighborhood of $f(t_2)$ with $\psi$. In this way, we
can eliminate the double points, separating via a fourth dimension. Thus we
obtain a ``resolution map''
\[
\Phi\co T^2= S_1\times S_2 \to \K_{4,1}.
\]

Finally, Budney showed that  $\Phi$ generates $H_{2}(\K_{4,1};\Z)$
and hence generates $\pi_{2}\K_{4,1}\cong\Z$ since
$\K_{4,1}$ is simply-connected \cite{BCSS}.

\subsection{Deform-spinning about a submanifold}\label{subsect:deformspin}
The spinning describes several methods of constructing higher-dimensional knots from
 lower-dimensional knots.
The most fundamental method, simply called \textit{spinning},
is due to Artin \cite{AR}.  It
has been generalized in various ways (a useful reference in this area is \cite{GF}).
In this paper, we will use the
\textit{deform-spinning about a submanifold}
(in \cite{GF} called \textit{frame deform-spinning}),
introduced in \cite{DR-spin}.

\begin{definition}\label{def:deformspin}
Suppose $M^k$ is a submanifold of $ S^{p}\subseteq S^{q}$ with
trivial tubular neighborhood $N\approx
M^k\times D^{p-k}\subseteq S^{p}$. Then $M$ has a trivial tubular
neighborhood $T \approx M^k\times D^{q-k}$ in $ S^q$. We can
write:
\[
(T,N) = M^k\times (D^{q-k}, D^{p-k}).
\]
Suppose $\Phi\co M^k\to\mathcal{K}_{q-k,p-k}$ is a smooth map. The
\textit{deform spun knot $\Sigma_\Phi$ about $M^k$} is the embedded
$p$-sphere in $S^q$, obtained from $S^p\subseteq S^q$ by replacing
the standard ball pair $\{x\}\times(D^{q-k}, D^{p-k})$ with
$\{x\}\times(D^{q-k},\Phi(x)( D^{p-k}))$ at each point $x\in
M^k$.
\end{definition}

This generalizes Litherland's \textit{deform-spinning} \cite{rick},
that corresponds to the case where $M^k$ is taken as $ S^1$
standardly embedded in $S^{p}\subseteq S^{p+2}$.

Deform-spinning about a submanifold is one of the most
generalized form of spinning.
For example, we can describe Artin's original
spinning as $\Sigma_\Phi$ where $M^k= S^1$ in $ S^2\subseteq S^3$
and $\Phi$ is a constant map. A \textit{super-spun} $p$-knot
\cite{CA} is given by deform-spinning about
$M^k= S^k\subseteq S^p$ with a constant map $\Phi$. Zeeman's
\textit{twist-spun} knot \cite{z2} is a deform-spun knot about
$M^k= S^1\subseteq S^2\subseteq S^3$ via the map
$\Phi\co S^1\to\mathcal{K}_{3,1}$; where $\Phi(\theta)$ is the
rotated image of a knotted arc about the $x$-axis by an angle of $\theta$.
Another variation is  Fox's \textit{roll-spinning} \cite{FX}.
Yet another extension
\textit{spinning of a knot about a projection of a knot}
uses a mapping of a manifold into $S^p$ \cite{DR-spin}.

\section{The main theorem}\label{sect:main}
Consider the $2$-torus $\T^2$ standardly embedded in $S^3$,
which is further included in $S^6$ in the standard manner:
$\T^2\subseteq S^3\subseteq S^6$.
At each point of $\T^2$ we
consider the normal $4$-disk to $\T^2\subseteq S^6$ and the
normal $1$-disk to $\T^2\subseteq S^3$, which form the standard
disk pair $(D^4,D^1)$.

We can deform-spin $ S^3\subseteq S^6$ about the torus
$\T^2\subseteq S^3\subseteq S^6$ 
with Budney's resolution map (\S\ref{subsubsect:resolution})
$\Phi\co\T^2\to\K_{4,1}$. In the normal plane at each point
$\tau\in\T^2$, we replace the standard disk pair $(D^4,D^1)$
with a new disk pair $(D^4,\Phi(\tau)(D^1))$.
We denote the resultant embedded $3$-sphere in
$S^6$ by $\Sigma_\Phi^3$. Namely, in $S^6$
\[
\Sigma_\Phi^3\quad=\quad
\overline{S^3\smallsetminus(\T^2\times D^1)}\quad\cup\quad\bigcup_{\tau\in \T^2}\Phi(\tau)(D^1).
\]

Let $J\co S^3\hookrightarrow S^6$ be a smooth embedding so that
$J( S^3)=\Sigma_\Phi^3\subseteq S^6$. Then, our main theorem is the
following:
\begin{thm}\label{thm:main}
$\Omega(J)=\pm1$; $J$ represents a generator of $C^3_3\cong\Z$.
\end{thm}

\section{A Seifert surface and the proof}\label{sect:seifert}
Consider the $2$-torus $\T^2\subseteq S^3\subseteq S^6$,
along which we performed the spinning, to be
\[
\T^2 = S^1_\theta\times S^1_\psi=
\{(\theta,\psi);\theta,\psi\in\R/2\pi\Z\}.
\]
At each point
$(\theta,\psi)\in \T^2$, we identify the normal $4$-disk
$D^4_{(\theta,\psi)}$ to $\T^2\subseteq S^6$ with the unit disk
$D^4$ in $\R^4=\{(x,y,z,w)\}$. Thus, $((\theta,\psi),(x,y,z,w))$
gives a coordinate system for a tubular neighborhood (diffeomorphic
to $\T^2 \times D^4$) of $\T^2 $ in $S^6$.

\begin{figure}[tbhp]
\begin{center}
\begin{picture}(80,80)
\put(17,57){$\mathbb{B}^2_{\theta,\psi}$}
\includegraphics[height=80mm,keepaspectratio]{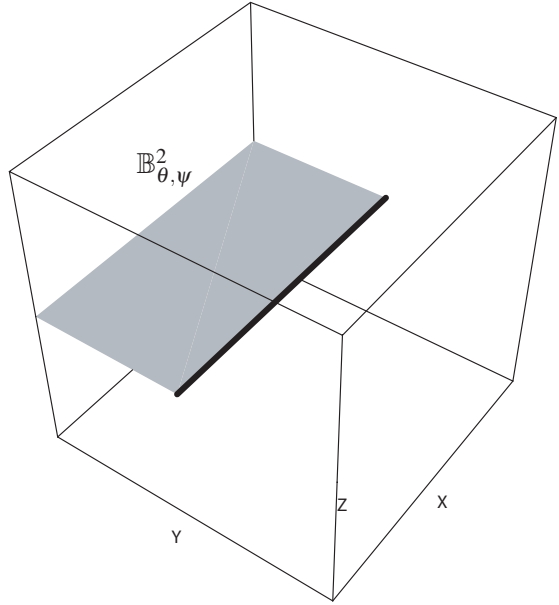}
\end{picture}
\end{center}
\caption{For any values of $(\theta,\psi)$, the standard Seifert surface of the unknot is as shown.
The center point of the cube has coordinates $(0,0,0,0)$;
the last coordinate $w$ is not shown.}\label{fig:trivial}
\end{figure}

Let $\mathbb{B}^4$ be the northern hemisphere of $S^4\subseteq S^6$,
which we think of the standard Seifert surface for the unknot $S^3\subseteq S^6$.
We can assume that
in each normal $4$-disk
$D^4_{(\theta,\psi)}$ to $\T^2 \subseteq S^6$,
the $4$-disk $\mathbb{B}^4$ is seen
as in Figure~\ref{fig:trivial},
which depicts the hyperplane section by $w=0$
of $D^4_{(\theta,\psi)}$ 
and where the intersection of
$\mathbb{B}^4$ and $D^4_{(\theta,\psi)}$ is shown in gray.
We denote this intersection $\mathbb{B}^4\cap D^4_{(\theta,\psi)}$
 by $\mathbb{B}^2_{\theta,\psi}$.

\begin{figure}[tbhp]
\begin{center}
\begin{picture}(80,80)
\put(15,56){$\puncT_{(0,0)}$}
\includegraphics[height=80mm,keepaspectratio]{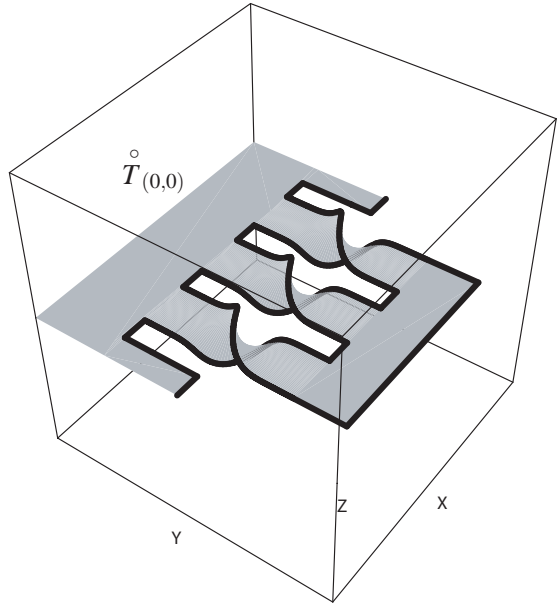}
\end{picture}
\end{center}
\caption{The Seifert surface for $\theta=\psi=0$ in $\{w=0\}$.}\label{fig:spun00}
\end{figure}

\subsection{A Seifert surface}\label{subsect:seifert}
\begin{figure}[tbhp]
\begin{center}
\includegraphics[width=80mm,keepaspectratio]{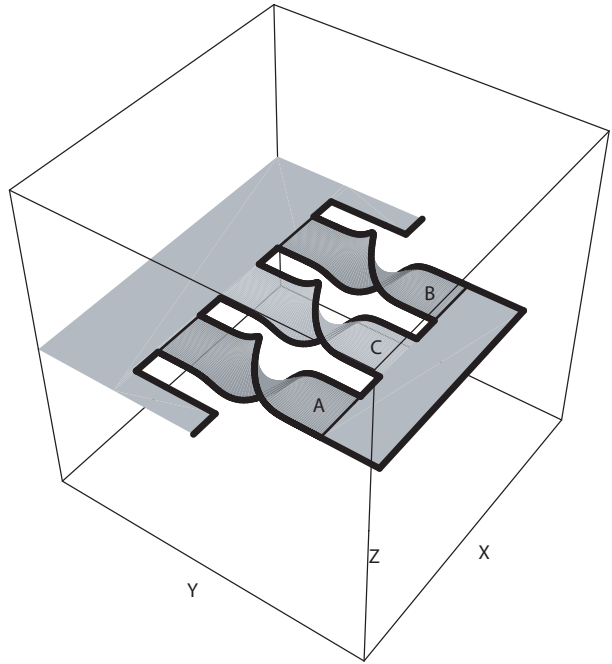}
\end{center}
\caption{The bands $A$ and $B$ of $\puncT_{(0,0)}$ are shown here in darkened gray,
the disk $C$ in light gray.}\label{fig:showBands}
\end{figure}

To construct a Seifert surface for our spun knot
$\Sigma_\Phi^3\subseteq S^6$, we first consider the punctured
$2$-dimensional torus in the $4$-disk $D^4_{(0,0)}$, as shown in
Figure~\ref{fig:spun00}. In Figure~\ref{fig:spun00}, $\Sigma_\Phi^3$
is the arc drawn with a thick line. We remark that although this arc
appears knotted since it is pictured in three-dimensional space, it is really unknotted in $D^4_{(0,0)}$. We view
$\puncT_{(0,0)}$ in a standard way as a $2$-disk $C$ with two
bands $A$ and $B$, as in Figure \ref{fig:showBands}.

In the normal $4$-disk $D^4_{(\theta,\psi)}$
for general $(\theta,\psi)\in \T^2$,
we consider the embedded punctured $2$-torus
$\puncT_{(\theta,\psi)}$,
defined as follows.
The punctured torus $\puncT_{(\theta,\psi)}$ coincides with
$\puncT_{(0,0)}$ on the disk $C$ and differs from $\puncT_{(0,0)}$
only on the two bands $A$ and $B$.
In the punctured torus $\puncT_{(\theta,\psi)}$,
the band $A$ has been replaced with $A_\theta$
and $B$ has been replaced with $B_\psi$:
$\puncT_{(\theta,\psi)} =  C \cup A_\theta \cup B_\psi$.
Now we describe the band $A_\theta$ in detail;
$B_\psi$ will be treated similarly.

\begin{figure}[tbhp]
\begin{center}
\includegraphics[width=80mm,keepaspectratio]{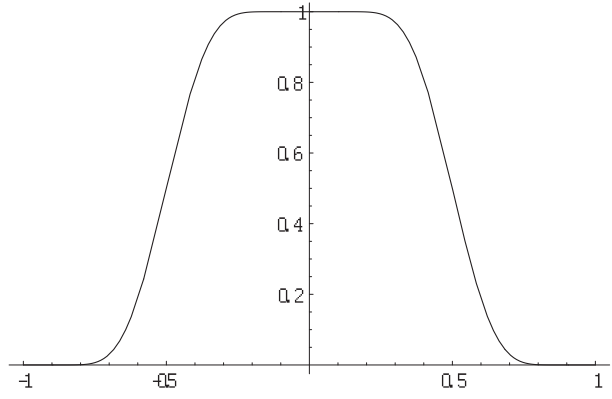}
\end{center}
\caption{The bump function}\label{fig:bump}
\end{figure}

We obtain $A_\theta$ by rotating $A$ around
the $2$-dimensional axis $\{w=z=0\}$
by an angle of $\gamma(y) \theta$, where
$\gamma(y)$ is a smooth bump function with $\gamma(y)=0$ in a
neighborhood of $\pm 1$ and $\gamma(y)=1$ for $-1/4\leq y\leq 1/4$
(see Figure~\ref{fig:bump}).
So $\gamma$ allows us to smoothly attach the band $A_\theta$ to $C$.
If we only consider the motion of
the edges of the band $A$ in the above process,
it corresponds to Budney's resolution process
(see \S\ref{subsubsect:resolution})
for one intersection point.

\begin{figure}[tbhp]
\begin{center}
\begin{picture}(80,80)
\put(14,56){$\puncT_{(\pi,0)}$}
\includegraphics[height=80mm,keepaspectratio]{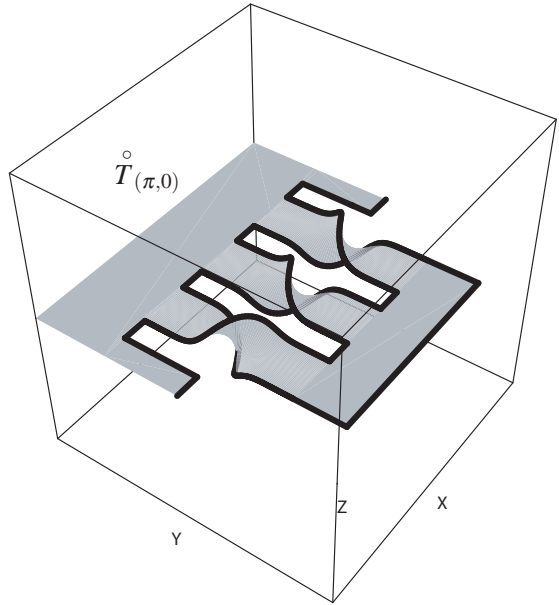}
\end{picture}
\end{center}
\caption{The Seifert surface for $\theta=\pi, \psi=0$ in $\{w=0\}$. 
Note the twist of the $A$ band is reversed if compared to Figure~\ref{fig:spun00}.}\label{fig:spun01}
\end{figure}

\begin{figure}[tbhp]
\begin{center}
\begin{picture}(80,80)
\put(13,56){$\puncT_{(\theta,\psi)}$}
\includegraphics[height=80mm,keepaspectratio]{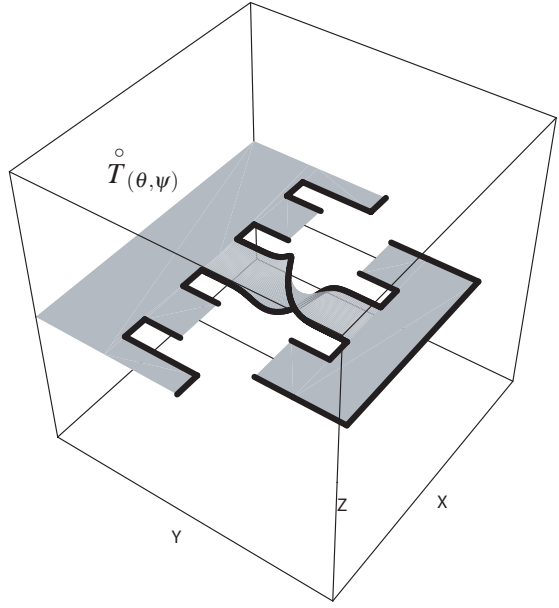}
\end{picture}
\end{center}
\caption{Generically, for  $\theta, \psi \neq 0, \mbox{ or } \pi$
the only parts of the two bands that appear in $\{w=0\}$ are the center lines.}\label{fig:spunGen}
\end{figure}

To clarify this construction, we show $\puncT_{(\pi,0)}$ in
Figure~\ref{fig:spun01}. In the generic case, the rotation will move
points not in the center of the band so that it has non-zero
$w$-coordinate and the points of $\puncT_{(\theta,\psi)}$ in
$\{w=0\}$ is shown in Figure~\ref{fig:spunGen}.

Putting together all the punctured tori
$\puncT_{(\theta,\psi)}\subseteq D^4_{(\theta,\psi)}$, we obtain the
embedded $4$-manifold
\[
W':=\displaystyle\bigcup_{(\theta,\psi)\in \T^2}\puncT_{(\theta,\psi)}
\quad\approx\T^2
\times\text{(the punctured torus)}\subseteq S^6.\]

Then,
\begin{eqnarray*}
W
&:=&\overline{\mathbb{B}^4\thinspace\smallsetminus\displaystyle\bigcup_{(\theta,\psi)\in \T^2}\thinspace\mathbb{B}^2_{(\theta,\psi)}}
\quad\cup\quad\displaystyle\bigcup_{(\theta,\psi)\in \T^2}\puncT_{(\theta,\psi)}\\
&=&\overline{\mathbb{B}^4\thinspace\smallsetminus\displaystyle\bigcup_{(\theta,\psi)\in \T^2}\thinspace\mathbb{B}^2_{(\theta,\psi)}}
\quad\cup\quad W'
\end{eqnarray*}
becomes a smoothly embedded $4$-manifold (a Seifert surface)
bounded by the knot
$\Sigma_\Phi^3\subseteq S^6$.
We abuse notation and identify
$W\subseteq S^6$ with the image of an embedding 
$\widetilde{J}\co W\hookrightarrow S^6$.

\subsection{The second homology group}\label{subsect:homology}
Now let us compute the second homology group
$H_2( W;\Z)$ of our Seifert surface $W$
and its intersection form.

\begin{figure}[tbhp]
\begin{center}
\includegraphics[width=80mm,keepaspectratio]{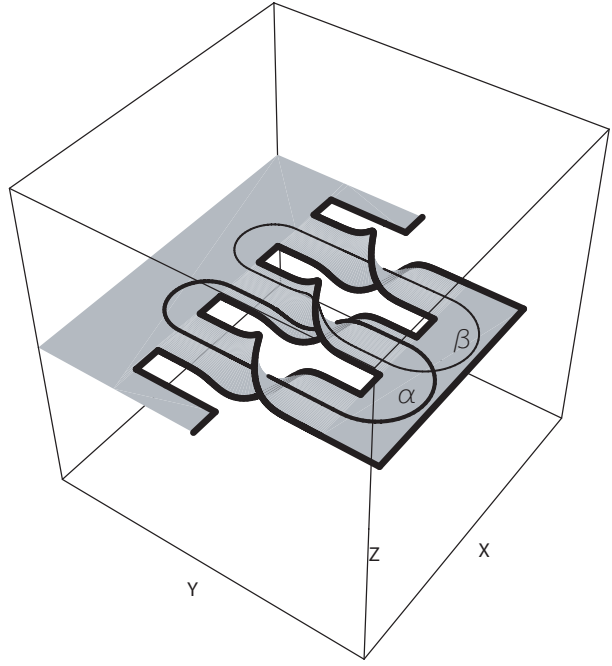}
\end{center}
\caption{The loops $\alpha$ and $\beta$.}\label{fig:ab}
\end{figure}

Since $W'$ is diffeomorphic to 
$\T^2\times\text{(the punctured $2$-torus)}$, we have
\[H_2(W')\approx\Z\oplus\Z\oplus\Z\oplus\Z\oplus\Z,\]
by the K\"uneth lemma. If we take two closed curves 
$\alpha$, $\beta$ 
in the punctured torus $\puncT_{(0,0)}$ at
$(0,0)\in S^1_\theta\times S^1_\psi$ as
shown in Figure~\ref{fig:ab}, then the above five $\Z$s are
generated by the five tori
$T_{\alpha\theta}:=\alpha\times S^1_\theta$,
$T_{\alpha\psi}:=\alpha\times S^1_\psi$,
$T_{\beta\theta}:=\beta\times S^1_\theta$,
$T_{\beta\psi}:=\beta\times S^1_\psi$ and
$\{\ast\}\times S^1_\theta\times S^1_\psi$, each embedded in $W'$.
Here $\ast$ is a point of $\puncT_{(0,0)}$.

Since $W$ is obtained by gluing 
$\overline{D^4\thinspace\smallsetminus\T^2\times D^2}$ ($\approx$ the $4$-disk) to $W'$
along a part of their boundaries, killing
$[\{\ast\}\times S^1_\theta\times S^1_\psi]\in H_2(W')$, we have
\begin{eqnarray*}
H_2(W)&\approx&\Z\oplus\Z\oplus\Z\oplus\Z\\
&\approx&
\langle[T_{\alpha\theta}]\rangle\oplus
\langle[T_{\beta\psi}]\rangle\oplus
\langle[T_{\alpha\psi}]\rangle\oplus
\langle[T_{\beta\theta}]\rangle,
\end{eqnarray*}
by the Mayer-Vietoris sequence.

Among the above representatives of $H_2(W)$,
the two pairs of ``complementary'' tori have intersection number $1$.
That is, the only non-zero intersection numbers are
\[[T_{\alpha\theta}]\cdot[T_{\beta\psi}]=[T_{\alpha\psi}]\cdot[T_{\beta\theta}]=1,\]
where $\cdot$ denotes the intersection pairing.
Hence, with respect to the above generators, the intersection form on $H_2(W)$
is expressed as
\[
\left(
  \begin{array}{cccc}
    0 & 1 & 0 & 0  \\
    1 & 0 & 0 & 0  \\
    0 & 0 & 0 & 1  \\
    0 & 0 & 1 & 0  \\
  \end{array}
\right).
\]
Note that the signature $\sigma(W)=0$.

\subsection{The normal Euler class}\label{subsect:euler}
We compute  the normal Euler class $e_{\widetilde{J}}$ for the
embedding $W\hookrightarrow S^6$ using an intersection argument.

Take a small generic perturbation $\hat{W}$ of $ W$ in $S^6$ and
put $F:=\hat{W}\cap W$. Since
$\Sigma_\Phi^3=\partial{W}\subseteq S^6$ has trivial normal bundle,
$F\subseteq\mathrm{Int}\thinspace W$ and the homology class $[F]\in
H_2(W)\approx H_2( W,\partial{W})$ is dual to the normal Euler class
$e_{\widetilde{J}}$.

With respect to the generators of $H_2(W)$ given in \S\ref{subsect:homology},
we represent the class $[F]$ as
\begin{eqnarray*}
[F]
&=&a_1[T_{\alpha\theta}]+a_2[T_{\beta\psi}]+a_3[T_{\alpha\psi}]+a_4[T_{\beta\theta}]\\
&=&(a_1,a_2,a_3,a_4)\in H_2(W)\cong\Z\oplus\Z\oplus\Z\oplus\Z.
\end{eqnarray*}
Then, with the intersection form in \S\ref{subsect:homology},
we have, for example
\begin{eqnarray*}
a_2=[F]\cdot[T_{\alpha\theta}].
\end{eqnarray*}
If we let $\hat{T}_{\alpha\theta}\subseteq\hat{W}$ be the perturbation of
$T_{\alpha\theta}$, then $a_2=[F]\cdot[T_{\alpha\theta}]$ is equal to the
intersection $F$ and $\hat{T}_{\alpha\theta}$ in $\hat{W}$. This is
further equal to the intersection of $W$ and $\hat{T}_{\alpha\theta}$,
since $F=\hat{W}\cap W$. Since $\hat{T}_{\alpha\theta}$ can be thought
of as a push-off of $T_{\alpha\theta}\subseteq W$ into
$S^6\smallsetminus W$,  we only need to  count the intersection of
$W$ and a push-off of $T_{\alpha\theta}$. By the same method, we can
compute $a_1$, $a_3$, $a_4$ and hence the class $[F]$ dual to the
normal Euler class $e_{\widetilde{J}}$.
We carry out this calculation below.

First, the punctured $2$-torus $\puncT_{(0,0)}$ at
$(\theta,\psi)=(0,0)$ lies in the $3$-dimensional hyperplane
$\{w=0\}$, the hyperplane section by $w=0$ of the normal $4$-disk
$ D^4_{(0,0)}$. Consider a vector field $\nu$ normal to the
punctured $2$-torus $\puncT_{(0,0)}$ in this $3$-dimensional hyperplane
$\{(x,y,z,0)\}$.
In each normal $4$-disk $D^4_{(\theta,0)}$, we consider
$(\theta,0)\times\alpha\subseteq\T^2\times D^4$
to be sitting in $\puncT_{(\theta,0)}$ and push it off
along the same normal
vector field $\nu$. This determines a push-off
$\hat{T}_{\alpha\theta}$ of the $2$-torus $T_{\alpha\theta}$ in $S^6$.

In each normal $4$-disk $D^4_{(\theta,0)}$, $W \cap D^4_{(\theta,0)}
= \puncT_{(\theta,0)}$. As we vary $\theta$, this punctured torus,
outside the $A$-band, lies in the $3$-dimensional hyperplane
$\{w=0\}$ and the $A$-band lies in
$\{(x,y,t\cos{\theta},t\sin{\theta})| t\in\R\}$. Therefore, the only
way that $W$ and $\hat{T}_{\alpha\theta}$ could intersect is when
$(\theta,\psi)=(\pi,0)$.

\begin{figure}[tbhp]
\begin{center}
\includegraphics[width=80mm,keepaspectratio]{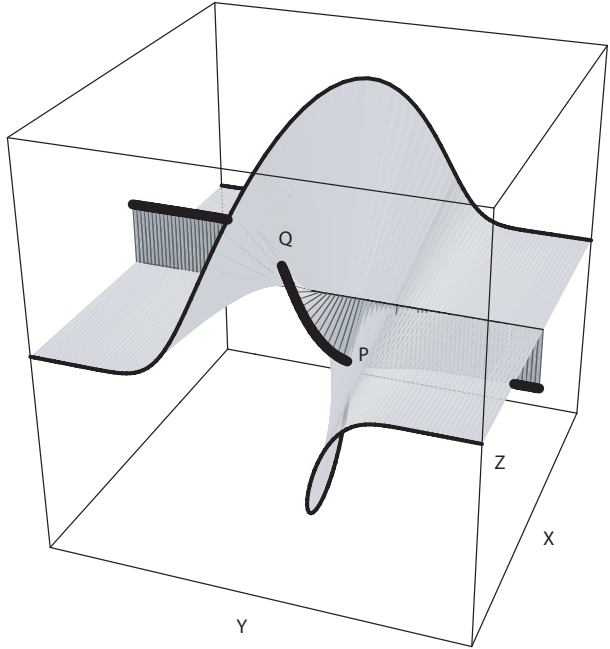}
\end{center}
\caption{Intersection of $\hat{T}_{\alpha\theta}$ and $W$.
We see the $A$-band of $\puncT_{(\pi,0)}\subseteq W$. 
A part of ${T}_{\alpha\theta}$ is shown as a straight thin line segment parallel to the $y$-axis.
The push-off $\hat{T}_{\alpha\theta}$ is a bold curve.
For visualization, we show the trace of this push-off.
}\label{fig:intersections}
\end{figure}

When $(\theta,\psi)=(\pi,0)$, in the normal $4$-plane
$D^4_{(\pi,0)}$, the whole $A$-band lies in the $3$-dimensional
hyperplane $\{w=0\}$. Figure~\ref{fig:intersections} depicts
the situation near the $A$-band in this $3$-dimensional plane, where
$\hat{T}_{\alpha\theta}$ is viewed as a ``half-twisted arc'' and the
$B$-band is an oppositely half-twisted band (Compare
Figures~\ref{fig:spun00} and \ref{fig:spun01}). We see that the
intersection of $W$ and $\hat{T}_{\alpha\theta}$ consists of the two
points $P$ and $Q$, seen in Figure~\ref{fig:intersections}.

\begin{figure}[tbhp]
\begin{center}
\includegraphics[width=\textwidth,keepaspectratio]{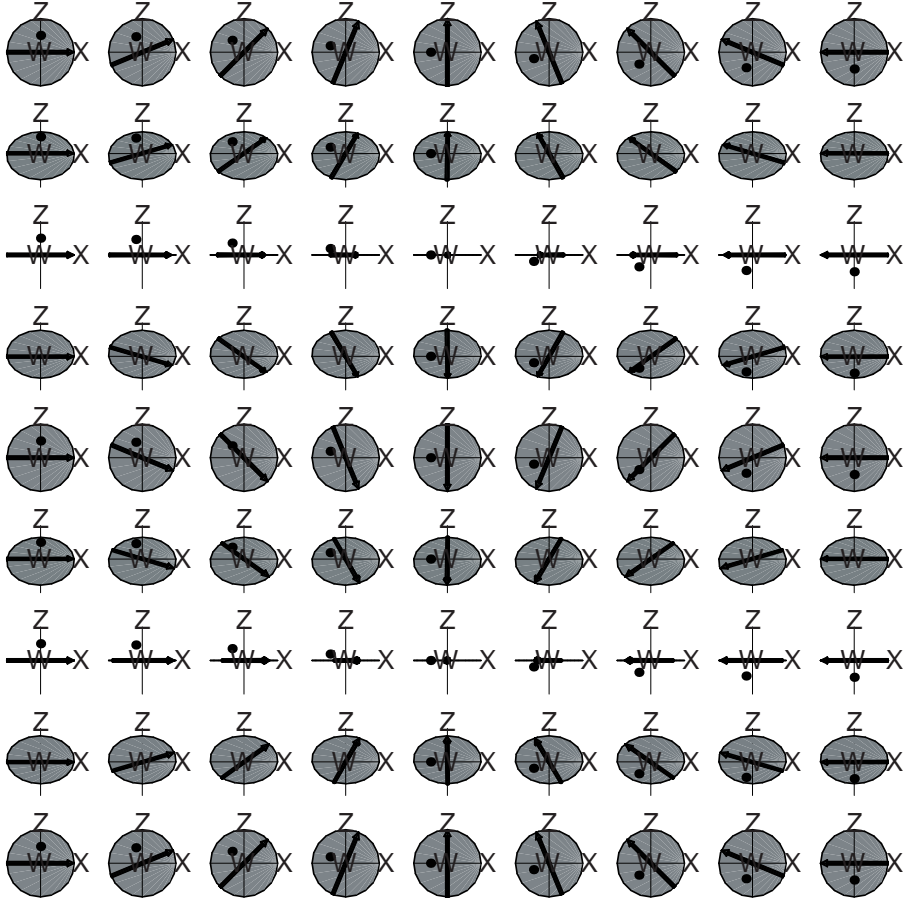}
\end{center}
\caption{Local pictures of $\hat{T}_{\alpha\theta}$ (dots) and $W$ (arrows).
The horizontal rows correspond to the $y$-direction and
the vertical columns correspond to the $\theta$-direction.
The disks in gray indicate the surfaces of revolution of the arrows
and only to guide visualization of these $3$-dimensional graphics.
They are not part of our construction.}\label{fig:movmov}
\end{figure}

However, it is not easy to see from Figure~\ref{fig:intersections}
the crucial fact that these intersections  actually have the same
sign. For this we use Figure~\ref{fig:movmov}.

Figure~\ref{fig:movmov} shows a $3$-dimensional sub-disk of $W$
intersecting a $2$-dimensional sub-disk of $\hat{T}_{\alpha\theta}$
transversely in $5$-dimensional space. The array of black dots
represent a square patch $S$ of $\hat{T}_{\alpha\theta}$. The arrow in
each figure represents a line segment  of an $A$-band in $W$; thus
the array of these segments represents a $3$-cube of $W$.

Not shown is a sixth coordinate --- the $\psi$ direction.  This is a
fourth coordinate for $W$ and the sixth coordinate of the ambient
space --- this coordinate will not be considered for analysis here.
What really matters is the normal bundle to $W$ and tangent plane
to the torus and these can be clearly understood using this figure.
The $\theta$ corresponds to twisting of the first band $A$ and is
independent from the twisting with respect to $\psi$.

In the fifth row of Figure~\ref{fig:movmov} the arrow and the dot for
each value of $y$ are both in the $xz$-plane in $\{w=0\}$. The
one-parameter family of these two-dimensional figures, when stacked,
give rise to the three dimensional Figure~\ref{fig:intersections}.

We will be concerned with the two  intersection points of our torus
$\hat{T}_{\alpha\theta}$ and the $4$-manifold $W$. In
Figure~\ref{fig:movmov} these points occur  where the black dot lies
exactly on the arrow
--- $P$ in row $5$ column $3$; and $Q$ in row $5$ column $7$. Choose an
orientation of the tangent bundle of $\hat{T}_{\alpha\theta}$ and
an orientation of the normal bundle of $W$.  Transversality assures
that, at each point of $W \cap \hat{T}_{\alpha\theta}$, we can
identify two-dimensional fibers of these bundles. The sign at the intersection
point $P$ is $+1$ if the two orientations agree, and $-1$ if not.

The square patch $S$ of $\hat{T}_{\alpha\theta}$ is flat in the 
$\theta xyz$-cube, thus its tangent plane coincides with this square $S$.
Specifically,
orient $\hat{T}_{\alpha\theta}$ with vectors $\tau_y$ in the $y$
direction (that is along the rows from left to right in
Figure~\ref{fig:movmov}) and $\tau_\theta$ the $\theta$ direction
(down the columns of Figure~\ref{fig:movmov}). This ordered pair
$(\tau_y,\tau_\theta)$ gives an orientation of the tangent plane of
$\hat{T}_{\alpha\theta}$. Next we chose a framing for the normal
bundle for W. The normal $2$-disks to $W$ are all represented as
disks normal to the arrow. In row $5$ column $1$ we orient the
normal disk $W$ by a vector pair $(\omega_1,\omega_2)$ where
$\omega_1$ is in the direction of the $z$-axis and $\omega_2$ in
direction of the $w$-axis. Continuing in the fifth row, we choose 
$\omega_1$ orthogonal to the arrow in the gray disk and $\omega_2$ in
the direction of the $w$-axis. (The details of the normal frames in
other rows will not need to be considered in our analysis.)

We first focus on $\tau_y$ at $P$ --- row $5$ column $3$. As we go
from left to right, the arrow rotates in the disk while the dot
rotates in the opposite direction. The dot
 from the positive side of the arrow (with respect to $\omega_1$)
to the negative side.  It follows that, at $P$, $\tau_y= -\omega_1$.
Similarly at $Q$ (row $5$ column $3$) the dot goes from the negative
side back to the positive side. Thus at $Q$, $\tau_y= +\omega_1$.
This can be seen also in Figure~\ref{fig:intersections}.

Next we consider $\tau_\theta$ at $P$. This second tangent direction
is downwards the third column. Note that as the dot passes $P$ it
goes from below the disk (that is, hidden from view) to above it
(visible).  This is in the direction of the positive $w$-axis. Thus
at $P$, $\tau_\theta=\omega_2$. At $Q$ we see the dot go from above
the disk to below it and so at $Q$, $\tau_\theta=-\omega_2$. This
information is not apparent in Figure~\ref{fig:intersections}.

In summary the intersection numbers at $P$ and $Q$ are both $-1$.
However, this sign depends on our choice of orientation of
$\hat{T}_{\alpha\theta}$, so we can only conclude that the signs at
$P$ and $Q$ are the same.

Thus, we conclude:
\begin{eqnarray*}
a_2=[F]\cdot[T_{\alpha\theta}]=\pm2.
\end{eqnarray*}
By the same argument for
$\beta$, $\psi$ and the $B$-band
instead of $\alpha$, $\theta$ and the $A$-band,
we have $a_1=\pm2$.

For the two tori $[T_{\alpha\psi}]$ and $[T_{\beta\theta}]$, we can
easily check that their push-offs along the same normal field $\nu$
 do not intersect
$W$ at all. Therefore, we have $a_3=a_4=0$.

Finally, the homology class dual to the desired normal Euler class is
\[
[F]=(\pm2,\pm2,0,0)\in H_2(W)\cong\Z\oplus\Z\oplus\Z\oplus\Z.
\]
Thus, we have

\begin{prop}\label{prop:euler}
$e_{\tilde{J}}=(\pm2,\pm2,0,0)\in H^2(W)$.
\end{prop}

\subsection{Proof of Theorem~\ref{thm:main}}\label{subsect:proof}
To prove Theorem~\ref{thm:main}, we have only to compute the
Haefliger invariant $\Omega(J)$ for our spun knot
$\Sigma_\Phi^3\subseteq S^6$
by using its Seifert surface $ W\subseteq S^6$ constructed in
\S\ref{subsect:seifert}.

By Theorem~\ref{thm:invariant},
together with Proposition~\ref{prop:euler}, we have
\begin{eqnarray*}
\Omega(J)
&=&-(\sigma(W)-e_{\widetilde{J}}\smile e_{\widetilde{J}})/8\\
&=&\pm(2\times2+2\times2)/8\\
&=&\pm1.
\end{eqnarray*}
This completes the proof of Theorem~\ref{thm:main}.\qed

\section{Remarks}\label{sect:remarks}
In view of the proof in \S\ref{subsect:proof},
we easily see the following.

\begin{remark}
In the construction of $\Sigma_\Phi^3\subseteq S^6$,
by using $\Phi'(\theta,\psi):=\Phi(m\theta,n\psi)$ ($m,n\in\Z$)
for spinning (i.e.\ if we change the speed of the resolutions),
we obtain the spun knot representing $mn$ times the generator
represented by $\Sigma_\Phi^3$ in $C^3_3$.
\end{remark}

All of our arguments are valid for higher-dimensional Haefliger
knots $C^{2k+1}_{4k-1}$, $k\ge2$. First of all, Budney's resolution
map is actually defined also in higher dimensions \cite[\S5]{budney}
and our construction of the spun knot works there. Furthermore,
since we also have higher-dimensional versions of
Theorem~\ref{thm:invariant} (see \cite[Theorem~2.3]{mT1} and
\cite[Theorem~5.1]{mT2}), the proof is directly extended in high
dimensions. Namely, we have: 

\begin{remark}
If we deform-spin $ S^{4k-1}\subseteq S^{6k}$ about
\[ S^{2k-1}\times S^{2k-1}\subseteq S^{4k-1}\subseteq S^{6k}\]
via the higher-dimensional Budney map \cite[\S5]{budney}
\[ S^{2k-1}\times S^{2k-1}\to\mathcal{K}_{2k+2,1}\]
(which corresponds to the case where we put $n=2k+2$ in Budney's
description of the generator of $\pi_{2n-6}\mathcal{K}_{n,1}$
\cite[\S5]{budney}), then the resultant spun knot represents a
generator of $C^{2k+1}_{4k-1}$ for $k\ge1$. All the steps of the
proof
 parallel  those of the
case $k=1$ (\S\ref{sect:seifert}) with little alteration.
\end{remark}


\end{document}